\documentclass[english]{article}
\usepackage[margin=2.5cm]{geometry}
\usepackage[hyphens]{url}
\usepackage{hyperref}
\usepackage[hyphenbreaks]{breakurl}
\usepackage{amsthm,amsfonts,amssymb,amsmath, color}
\usepackage{graphicx}
\usepackage{bbm}
\usepackage{tikz}

\newtheorem{theorem}{Theorem}[section]
\newtheorem{proposition}[theorem]{Proposition}
\newtheorem{lemma}[theorem]{Lemma}

\newtheorem{fact}[theorem]{Fact}

\newtheorem{corollary}[theorem]{Corollary}
\newtheorem{definition}[theorem]{Definition}

\theoremstyle{definition}

\numberwithin{equation}{section}
\numberwithin{figure}{section}
\allowdisplaybreaks

\usepackage{babel}
\usepackage[T1]{fontenc}
\usepackage[utf8]{inputenc}

\newcommand{\su}{\subseteq}

\newcommand{\Fp}{\mathbb{F}_p}
\newcommand{\Fpn}{\mathbb{F}_p^n}
\newcommand{\Z}{\mathbb{Z}}

\begin{document}
\title{New lower bounds for three-term progression free sets in $\mathbb{F}_p^n$}
\author{Christian Elsholtz\thanks{Institute for Analysis and Number Theory, TU Graz, Austria. Email: \url{elsholtz@math.tugraz.at}.
Supported by the joint FWF-ANR project Arithrand: FWF: I 4945-N and ANR-20-CE91-0006, and FWF grant DOC 183.}\and Laura Proske\thanks{TU Graz, Austria. Email: \url{proske@student.tugraz.at}.}\and Lisa Sauermann\thanks{Institute for Applied Mathematics, University of Bonn, Germany. Email: \url{sauermann@iam.uni-bonn.de}. Supported by the DFG Heisenberg Program.}}

\maketitle

\begin{abstract}\noindent
We prove new lower bounds on the maximum size of sets $A\subseteq \mathbb{F}_p^n$ or $A\subseteq \Z_m^n$ not containing three-term arithmetic progressions (consisting of three distinct points). More specifically, we prove that for any fixed integer $m\ge 2$ and sufficiently large $n$ (in terms of $m$), there exists a three-term progression free subset $A\subseteq \mathbb{Z}_m^n$ of size $|A|\ge (cm)^n$ for some absolute constant $c>1/2$. Such a bound for $c=1/2$ can be obtained with a classical construction of Salem and Spencer from 1942, and improving upon this value of $1/2$ has been a well-known open problem (our proof gives $c= 0.54$).

Our construction relies on finding a subset $S\subseteq \mathbb{Z}_m^2$ of size at least $(7/24)m^2$ with a certain type of reducibility property. This property allows us to ``lift'' $S$ to a three-term progression free subset of $\mathbb{Z}_m^n$ for large $n$ (even though the original set $S\subseteq \mathbb{Z}_m^2$ does contain three-term arithmetic progressions).
\end{abstract}

\section{Introduction}

Studying the maximum possible size of subsets of $\{1,\dots,N\}$ or of $\mathbb{F}_p^n$ without three-term arithmetic progressions is one of the most central problems in additive combinatorics. In a breakthrough result in 2017, Ellenberg and Gijswijt \cite{ellenberg-gijswijt} proved that for any prime $p\ge 3$, every subset $A\su \Fpn$ not containing a three-term arithmetic progression (i.e.\ not containing three distinct elements $x,y,z$ with $x-2y+z=0$) has size at most $(c_pp)^n$ for some constant $c_p<1$ (where $c_p$ converges to $0.841\dots$ for $p\to \infty$, see \cite[Eq.\ (4.11)]{blasiak-et-al}). In the integer setting, in a more recent breakthrough Kelley and Meka \cite{kelley-meka} proved that for $N\ge 3$, every subset $A\su \{1,\dots,N\}$ not containing a three-term arithmetic progression has size at most $N\exp(-c(\log N)^{1/11})$ for some absolute constant $c>0$ (afterwards, this bound was improved to $N\exp(-c(\log N)^{1/9})$ using a modification of their method by Bloom and Sisask \cite{bloom-sisask}). This matches the shape of a classical lower bound for this problem due to Behrend \cite{behrend}, which is of the form $N\exp(-C(\log N)^{1/2})$ for some absolute constant $C$. See also \cite{peluse} for a recent survey on this problem (both in $\{1,\dots,N\}$ and in $\mathbb{F}_p^n$), and its analogues for longer arithmetic progressions and polynomial progressions.

To obtain lower bonds for the maximum size of three-term progression free subsets of $\mathbb{F}_p^n$ for large $n$, one can adapt a classical construction of Salem and Spencer \cite{salem-spencer} from 1942, see \cite[Theorem 2.13]{alon-shpilka-umans}. Alternatively, one can also adapt Behrend's construction \cite{behrend} from 1946 to this setting, as noted by Tao and Vu in their book on additive combinatorics \cite[Exercise 10.1.3]{tao-vu} and also observed by Alon (see \cite[Lemma 17]{fox-pham}). Both of these constructions show that for every fixed prime $p$, there exists a subset $A\su \Fpn$ of size $|A|\ge ((p+1)/2)^{n-o(n)}$ without a three-term arithmetic progression (here, the asymptotic notation $o(n)$ is for $n\to \infty$ with $p$ fixed). The best quantitative bound for the $o(n)$-term in this statement is due to the first author and Pach \cite[Theorem 3.10]{elsholtz-pach}. For small primes, better bounds were obtained with specific constructions depending on the particular prime. For example, in $\mathbb{F}_3^n$ a lower bound of $2.2202^{n-o(n)}$ was recently obtained by Romera--Paredes et al.\ \cite{google-deep-mind} using artifical intelligence building upon traditional methods from previous bound \cite{calderbank-fishburn,edel,tyrrell}. Naslund \cite{naslund} informed us about forthcoming work with an approach related to Shannon capacity, proving a lower bound of $2.2208^{n-o(n)}$ for the maximum size of three-term progression free subsets of $\mathbb{F}_3^n$. In $\mathbb{F}_5^n$, the best known lower bounds is $(35^{1/3})^{n-o(n)}$ due to the first two authors, Pollak, Lipnik and Siebenhofer (see \cite{elsholtz-lipnik-siebenhofer}), note that $35^{1/3}\approx 3.271$.

Comparing the upper and lower bounds for the maximum possible size of a three-term progression free subset of $\mathbb{F}_p^n$ for a fixed (reasonably large) prime $p$ and large $n$, there is still a large gap. Both of these bounds are roughly of the form $(cp)^n$ with $0<c<1$, but with different values of $c$. For the upper bound, the best known constant due to Ellenberg--Gijswijt \cite{ellenberg-gijswijt} is $c\approx 0.85$ (when the fixed prime is large enough), whereas for the lower bound the constant $c=1/2$ from the Salem--Spencer construction \cite{salem-spencer} (or alternatively Behrend's construction \cite{behrend}) has not been improved in more than 80 years. In this paper, we finally improve this constant in the lower bound to be strictly larger than $1/2$.

\begin{theorem}\label{thm-main}
There is a constant $c>1/2$ such that for every prime $p$ and every sufficiently large positive integer $n$ (sufficiently large in terms of $p$), there exists a subset $A\su \Fpn$ of size $|A|\ge (cp)^n$ not containing a three-term arithmetic progression.
\end{theorem}

Our proof shows that one can take any $c< \sqrt{7/24}$, for example $c=0.54$. While quantitatively, this may not seem like a large improvement over $1/2$, it is the first qualitative improvement over the constant $1/2$ from the constructions of Salem--Spencer and Behrend from the 1940's. Both of these constructions lead to three-term progression free subsets of $\Fpn$ only consisting of vectors with all entries in $\{0,1,\dots,(p-1)/2\}$, i.e.\ they only use roughly half of the available elements in $\Fp$ in each coordinate. The restriction of all entries to $\{0,1,\dots,(p-1)/2\}$ is crucial in these constructions, as it ensures that there is no ``wrap-around'' over $\Fp$. However, such an approach cannot be used to obtain three-term progression free subset $A\su \Fpn$ of size $|A|> ((p+1)/2)^n$, so $c=1/2$ is a significant barrier for this problem (and it may be considered a surprise that it is actually possible to obtain a constant $c>1/2$).

In fact, our construction proving Theorem \ref{thm-main} works more generally in $\Z_m^n$ for any integer $m\ge 2$ (not necessarily prime) and sufficiently large $n$, again showing that there exists a three-term progression free subset of $\Z_m^n$ of size at least $(cm)^n$. We recall that a three-term arithmetic progression in $\Z_m^n$ consists of three distinct elements $x,y,z\in \Z_m^n$ with $x-2y+z=0$ (see e.g.\ the breakthrough work of Croot--Lev--Pach \cite{croot-lev-pach} proving upper bounds for three-term progression free subsets of $\Z_4^n$). For odd $m$, the best previous lower bound for this problem has also been of the form $((m+1)/2)^{n-o(n)}$ (based on the Salem--Spencer construction \cite{salem-spencer} or the Behrend construction \cite{behrend}), and for even $m$ the best previous lower bound has been $((m+2)/2)^{n-o(n)}$ due to the first author and Pach \cite[Theorem 3.11]{elsholtz-pach}.

\begin{theorem}\label{thm-main-2}
There is a constant $c>1/2$ such that for every integer $m\ge 2$ and every sufficiently large positive integer $n$ (sufficiently large in terms of $m$), there exists a subset $A\su \Z_m^n$ of size $|A|\ge (cm)^n$ not containing a three-term arithmetic progression.
\end{theorem}

Note that Theorem \ref{thm-main-2} immediately implies Theorem \ref{thm-main}. So it suffices to prove Theorem \ref{thm-main-2}.

Our construction of the set $A\su \Z_m^n$ in Theorem \ref{thm-main-2} is, in some sense, a higher-dimensional version of the Salem--Spencer construction \cite{salem-spencer}. Generalizing the argument of Salem and Spencer requires a certain reducibility notion for subsets of $\Z_m^d$. In the one-dimensional case (i.e.\ for $d=1$), a similar type of reducibility argument was used by the first author and Lipnik \cite{elsholtz-lipnik} to construct subsets of $\Fpn$ without three points on a common line for $p\in \{11, 17, 23, 29, 41\}$ and large $n$, and also by the first author, Klahn and Lipnik \cite{elsholtz-klahn-lipnik}
to construct subsets of $\Z_m^n$ avoiding arithmetic progressions of a given length $k> 3$ (other constructions for this problem with different methods were given in \cite{elsholtz-et-al-line-free,fuehrer,lin-wolf}). A higher-dimensional version of this type of argument first appeared in the above-mentioned work of the first two authors, Pollak, Lipnik and Siebenhofer (see \cite{elsholtz-lipnik-siebenhofer}) constructing three-term progression free subsets of $\mathbb{F}_5^n$. However, the crucial difficulty for implementing this strategy for larger  $m$ is to find suitable higher-dimensional building blocks for this construction, i.e.\ suitable large subsets of $\Z_m^d$ with the required reducibility condition. In that work, a suitable subset of $\mathbb{Z}_5^3$ was found via a computer search. While for any given $m$ and $d$ one can try to find suitable subsets  of $\Z_m^d$ via computer searching, the difficulty in proving Theorems \ref{thm-main} and \ref{thm-main-2} is establishing the existence of such subsets for all $m$.

The problem of finding the maximum possible size of a three-term progression free subset has also been studied in the context of general finite abelian groups.
Brown and Buhler \cite{brown-buhler} proved in 1982 that every three-term progression free subset of an abelian group $G$ has size $o(|G|)$. Frankl, Graham and R\"odl \cite{frankl-graham-rodl} gave a different proof of this statement in 1987, and further upper bounds were obtained by Meshulam \cite{meshulam} and Lev \cite{lev:2004}. Every finite abelian group $G$ can be written as a product of factors of the form $\Z_m^n$, and by a simple product-construction the maximum possible size of a three-term progression free subset of $G$ is at least the product of the corresponding maximum subset sizes for these factors. Thus, as long as there are factors $\Z_m^n$ with $n$ sufficiently large in terms of $m$, Theorem \ref{thm-main-2} can also be used to obtain lower bounds for more general finite abelian groups.

We will give the proof of Theorem \ref{thm-main-2} in Section \ref{sect-overview}, apart from postponing the crucial construction of the above-mentioned building blocks to Section \ref{sect-constr-non-midpoint-deg-set} (recall that Theorem \ref{thm-main} follows from Theorem \ref{thm-main-2}). We end with some concluding remarks in Section \ref{sect-concluding}.

\textit{Acknowledgements.} We would like to thank Eric Naslund for interesting discussions.

\section{Proof Overview}
\label{sect-overview}

As mentioned in the introduction, our construction can be viewed as a higher-dimensional generalization of the classical Salem--Spencer construction. However, the higher-dimensionality makes the construction more intricate, and one of the main difficulties is finding higher-dimensional building blocks of that lead to a better bound than the classical Salem--Spencer construction with one-dimensional building blocks. The following definitions are used to describe the relevant properties of these building blocks of our construction.

\begin{definition}\label{defi-non-mid-point}
For a subset $S\su \Z_m^d$, a point $y\in S$ is called a \emph{non-mid-point} of $S$ if there is no solution to the equation $2y=x+z$ with distinct $x,z\in S$. 
\end{definition}

Note that for odd $m$ and a subset $S\su \Z_m^d$, a point $y\in S$ is a non-mid-point of $S$ if and only if the only solution to $2y=x+z$ with $x,z\in S$ is $x=z=y$. For even $m$, on the other hand, if $y\in S$ is a non-mid-point of $S$, there can be additional solutions to $2y=x+z$ with $x,z\in S$ besides $x=z=y$. To describe these, we introduce the following definition.

\begin{definition}
For even $m$, we say that two points $x,y\in \Z_m^d$ are \emph{cousins} if $2x=2y$.
\end{definition}

Note that every point in $\Z_m^d$ (for even $m$) is a cousin of itself. Furthermore, note that the relation of being cousins is an equivalence relation on $\Z_m^d$ (in particular, it is transitive), so it divides $\Z_m^d$ into equivalence classes. Also note that for a subset $S\su \Z_m^d$ and two cousins $x,y\in S$, the point $x$ is a non-mid-point of $S$ if and only if the point $y$ is a non-mid-point of $S$.

Finally, note that for even $m$ and a subset $S\su \Z_m^d$, a point $y\in S$ is a non-mid-point of $S$ if and only if all the solutions to $2y=x+z$ with $x,z\in S$ are such that $x=z$ is a cousin of $y$.

\begin{definition}\label{defi-non-mid-point-reducible}
A subset $S\su \Z_m^d$ is called \emph{non-mid-point reducible} if for every non-empty subset $S'\su S$, some point $y\in S'$ is a non-mid-point of $S'$.
\end{definition}

A subset $S\su \Z_m^d$ is non-mid-point reducible, if and only if it can be reduced to the empty set by repeatedly removing non-mid-points. More precisely, $S$ is non-mid-point reducible, if and only if some point $y_1\in S$ is a non-mid-point of $S$, and some point $y_2\in S\setminus \{y_1\}$ is a non-mid-point of $S\setminus \{y_1\}$, and some point $y_3\in S\setminus \{y_1,y_2\}$ is a non-mid-point of $S\setminus \{y_1,y_2\}$, and so on, until the final set $S\setminus \{y_1,y_2,\dots,y_{|S|}\}$ is empty. 

One can use non-mid-point reducible subsets of $\Z_m^d$ (for some fixed $m$ and $d$) to construct three-term progression free subsets of $\Z_m^n$ (for large $n$). Indeed, by appropriately generalizing the arguments of Salem and Spencer, one can show the following proposition.

\begin{proposition}\label{proposition-1}
Let $m$ and $d$ be fixed positive integers, and let $S\su \Z_m^d$ be a non-mid-point reducible subset of $\Z_m^d$. Then, there exists a subset $A\su \Z_m^n$ of size $|A|\ge |S|^{n/d-O(\log n)}$ not containing a three-term arithmetic progression.
\end{proposition}

Here, the asymptotic notation $O(\log n)$ is for $n\to \infty$ with $m$ and $d$ fixed. The Salem--Spencer construction for three-term progression free subsets of $\Fpn$ can be viewed as applying this proposition to $m=p$ and $d=1$ with $S=\{0,1,\dots,(p-1)/2\}\su \Z_p$. The above-mentioned work \cite{elsholtz-lipnik-siebenhofer} on three-term progression free subsets of $\mathbb{F}_5^n$ used the argument underlying this proposition for $d=3$. 

We will prove Proposition \ref{proposition-1} at the end of this section. In order to prove Theorem \ref{thm-main} using this proposition, we of course also need to show that there exist reasonably large non-mid-point reducible subsets $S\su \Z_m^d$ for some $d$. This is the content of the following proposition, and is one of the main difficulties of this paper.

\begin{proposition}\label{proposition-2}
For any positive integer $m$, there exists a non-mid-point reducible subset $S\su \Z_m^2$ of size $|S|\ge (7/24)\cdot m^2$.
\end{proposition}

We prove Proposition \ref{proposition-2} in Section \ref{sect-constr-non-midpoint-deg-set}. We remark that it is fairly easy to find a non-mid-point reducible subset $S\su \Z_m^2$ of size $|S|\ge \lfloor (m+2)/2\rfloor^2$ (for example, one can take $S=\{0,1,\dots,\lfloor m/2\rfloor\}^2$). However, this would only lead to a bound of $|A|\ge \lfloor (m+2)/2\rfloor^{n-O(\log n)}$ in Theorem \ref{thm-main-2} (compare \cite[Theorem 3.11]{elsholtz-pach}), i.e.\ it would not give a constant $c>1/2$. The construction of the set $S\su \Z_m^2$ of size $|S|\ge (7/24)\cdot m^2$ in Proposition \ref{proposition-2} is much more involved, and it also takes more work to show that this set $S$ is indeed non-mid-point reducible.

Combining Propositions \ref{proposition-1} and \ref{proposition-2}, we can now show that Theorem \ref{thm-main-2} holds for any constant $c<\sqrt{7/24}$ (note that $\sqrt{7/24}\approx 0.54006...$). In fact, for any fixed $m\ge 2$, we obtain a three-term progression free subsets $A\su \mathbb{Z}_m^n$ of size $|A|\geq (\sqrt{7/24}\cdot m)^{n-O(\log n)}$.

\begin{proof}[Proof of Theorem \ref{thm-main-2}, assuming Propositions \ref{proposition-1} and \ref{proposition-2}]

Noting that $\sqrt{7/24}>\sqrt{1/4}=1/2$, choose a constant $c$ with  $1/2<c<\sqrt{7/24}$. By Proposition \ref{proposition-2}, for every positive integer $m$ there exists a non-mid-point reducible subset $S\su \Z_m^2$ of size $|S|\ge (7/24)\cdot m^2$. Now, by Proposition \ref{proposition-1}, for large $n$ there is a progression-free subset $A\su \Z_m^n$ of size
\[|A|\ge |S|^{n/2-O(\log n)}=\left(|S|^{1/2}\right)^{n-O(\log n)}\ge \left(\sqrt{7/24} \cdot m\right)^{n-O(\log n)}\ge \left(\sqrt{7/24}  \cdot m\right)^{(1-o(1))n},\]
where the $o(1)$-term converges to zero as $n\to \infty$ (while $m$ is fixed). Given that $\sqrt{7/24} \cdot m>cm$, for all sufficiently large $n$ (in terms of $m$ and the fixed constant $c$), we can therefore find a progression-free subset $A\su \Z_m^n$ of size $|A|\ge (cm)^n$.
\end{proof}

We end this section by proving Proposition \ref{proposition-1}. First, we prove the proposition in the case that $m$ is odd (since this case involves fewer technicalities).

\begin{proof}[Proof of Proposition \ref{proposition-1} for odd $m$]
First, note that the statement is trivially true if $|S|\le 1$, so we may assume that $|S|\ge 2$.

Furthermore, noting that $|S|\le m^d$, we may assume that $n$ is divisible by $d\cdot |S|$. Indeed, assuming that we have proved the desired statement under this assumption, for any large positive integer $n$ (large in terms of $m$ and $d$) we can take some positive integer $n'$ with $n-d\cdot |S|\le n'\le n$ such that $n'$ is divisible by $d\cdot |S|$, and find a three-term progression free subset $A\su \Z_m^{n'}$ of size $|A|\ge |S|^{n'/d-O(\log n')}\ge |S|^{n/d-|S|-O(\log n)}=|S|^{n/d-O(\log n)}$. Embedding $\Z_m^{n'}$ into $\Z_m^{n}$ by appending $n-n'$ zeroes to each point, this also gives a three-term progression free subset of $\Z_m^{n}$ of the same size.

So let us from now on assume that $n=d\cdot |S|\cdot \ell$ for some large positive integer $\ell$. We may now view $\Z_m^{n}$ as a product $\Z_m^{n}=\Z_m^{d}\times \dots \times \Z_m^{d}$ of $\ell |S|$ factors $\Z_m^{d}$. We define
\[A\su S\times \dots \times S\su \Z_m^{d}\times \dots \times \Z_m^{d}=\Z_m^{n}\]
to be the set of all $(\ell |S|)$-tuples $(w_1,\dots,w_{\ell |S|})\in S\times \dots \times S=S^{\ell |S|}$ such that each point $w\in S$ appears exactly $\ell$ times among $w_1,\dots,w_{\ell |S|}$. Then we have
\[|A|=\binom{\ell |S|}{\ell,\ell,\dots,\ell}\ge \frac{|S|^{\ell |S|}}{(\ell |S|+1)^{|S|}}\ge \frac{|S|^{n/d}}{(n+1)^{|S|}}\ge |S|^{n/d-|S|\log_2 (n+1))}\ge |S|^{n/d-O(\log n)},\]
using that the multinomial coefficient in the second term is the largest of the at most $(\ell |S|+1)^{|S|}$ multinomial coefficients with the same value $\ell |S|$ at the top and $|S|$ parts at the bottom\footnote{A more careful analysis of the multinomial coefficient with Stirling's formula $\sqrt{2\pi}\cdot k^{k+1/2}e^{-k} \leq k! \leq e \cdot k^{k+1/2}e^{-k} $ (for all $k \geq 1$) yields $|A|\geq |S|^{n/d}/(e^2 \ell)^{(|S|-1)/2}$.}, and using that $|S|\le m^d$ with $m$ and $d$ being fixed for the asymptotic term $O(\log n)$.

It remains to prove that the set $A$ does not contain a three-term arithmetic progression. Suppose for contradiction that $x=(x_1,\dots,x_{\ell |S|})$ and $y=(y_1,\dots,y_{\ell |S|})$ and $z=(z_1,\dots,z_{\ell |S|})$ are three distinct elements of $A$ with $x-2y+z=0$. Then we have $2y_i=x_i+z_i$ for all $i=1,\dots,\ell |S|$.

For an index $i\in \{1,\dots,\ell |S|\}$, let us say that $i$ is \emph{boring} if $x_i=y_i=z_i$, and that $i$ is \emph{interesting} otherwise. As $x,y,z$ are distinct, not all indices $i\in \{1,\dots,\ell |S|\}$ can be boring. Let
\[S'=\{y_i\mid i\in \{1,\dots,\ell |S|\}\text{ interesting}\}\su S,\]
then $S'\ne\emptyset$, since there must be at least one interesting  index $i\in \{1,\dots,\ell |S|\}$.

By definition of $S'$, for every interesting index $i\in \{1,\dots,\ell |S|\}$ we have $y_i\in S'$. We claim that furthermore, for every interesting index $i\in \{1,\dots,\ell |S|\}$, we also have $x_i\in S'$. Indeed, if $x_i=y_i$, then we trivially have $x_i=y_i\in S'$. Otherwise, if $x_i\ne y_i$, then there are $\ell$ indices $j\in \{1,\dots,\ell |S|\}\setminus \{i\}$ with $y_j=x_i$ (as $x_i\in S$ appears exactly $\ell$ times among $y_1,\dots,y_{\ell |S|}$). At most $\ell-1$ of these indices $j$ can have the property that $x_j=x_i$, so there must be an index $j\in \{1,\dots,\ell |S|\}\setminus \{i\}$ such that $y_j=x_i$ and $x_j\ne x_i$. Then we must have $x_j\ne y_j$, so the index $j$ is also interesting and hence $x_i=y_j\in S'$ as claimed.

This shows that for every interesting index $i\in \{1,\dots,\ell |S|\}$, we have $x_i\in S'$. Analogously, one can also show $z_i\in S'$ for every interesting index $i\in \{1,\dots,\ell |S|\}$.

As the set $S$ is mid-point-reducible, the non-empty subset $S'\su S$ must contain some point that is a non-mid-point of $S'$. So there must be some interesting index $i\in \{1,\dots,\ell |S|\}$, such that $y_i$ is a non-mid-point of $S'$. On the other hand, we have $2y_i=x_i+z_i$ and $x_i,z_i\in S'$, so by Definition \ref{defi-non-mid-point} this means that $x_i=z_i$. Hence $2y_i=x_i+z_i=2x_i$, and as $m$ is odd, this implies $y_i=x_i$. So we have $x_i=y_i=z_i$, which is a contradiction to $i$ being an interesting index. This contradiction shows that the set $A$ is indeed three-term progression free.
\end{proof}

Now, we prove Proposition \ref{proposition-1} for even $m$. The proof is essentially the same as above, but one needs to be more careful by handling cousins in a suitable way (since these can lead to additional solutions to $2y=x+z$ for a  non-mid-point $y$ besides $x=z=y$).

\begin{proof}[Proof of Proposition \ref{proposition-1} for even $m$]
For the same reasons as in the previous proof for odd $m$, we may assume that $|S|\ge 2$ and that $n$ is divisible by $d\cdot |S|$. So let $n=d\cdot |S|\cdot \ell$ for some large positive integer $\ell$. We again view $\Z_m^{n}$ as a product $\Z_m^{n}=\Z_m^{d}\times \dots \times \Z_m^{d}$ of $\ell |S|$ factors $\Z_m^{d}$, and define
\[A\su S\times \dots \times S\su \Z_m^{d}\times \dots \times \Z_m^{d}=\Z_m^{d}\]
to be the set of all $(\ell |S|)$-tuples $(w_1,\dots,w_{\ell |S|})\in S\times \dots \times S=S^{\ell |S|}$ such that each point $w\in S$ appears exactly $\ell$ times among $w_1,\dots,w_{\ell |S|}$. Then we again have $|A|\ge |S|^{n/d-O(\log n)}$ by the same calculation as in the previous proof for odd $m$.

Again, it remains to prove that the set $A$ does not contain a three-term arithmetic progression. Suppose for contradiction that $x=(x_1,\dots,x_{\ell |S|})$ and $y=(y_1,\dots,y_{\ell |S|})$ and $z=(z_1,\dots,z_{\ell |S|})$ are three distinct elements of $A$ with $x-2y+z=0$, then $2y_i=x_i+z_i$ for $i=1,\dots,\ell |S|$.

For an index $i\in \{1,\dots,\ell |S|\}$, let us say that $i$ is \emph{boring}, if $x_i=z_i$ and $y_i$ is a cousin of $x_i=z_i$, and let us say that $i$ is \emph{interesting} otherwise. Again, not all indices $i\in \{1,\dots,\ell |S|\}$ can be boring, since then we would in particular have $x_i=z_i$ for all $i\in \{1,\dots,\ell |S|\}$, meaning that $x=z$, which contradicts $x,y,z$ being distinct. So there must be at least one interesting index $i\in \{1,\dots,\ell |S|\}$. Let
\[S'=\{w\in S\mid w\text{ is a cousin of }y_i\text{ for some interesting }i\in \{1,\dots,\ell |S|\}\},\]
then $S'\su S$ and $S'\ne\emptyset$.

Clearly, for every interesting index $i\in \{1,\dots,\ell |S|\}$, we have $y_i\in S'$ (recalling that $y_i$ is a cousin of itself). We again claim that we also have $x_i\in S'$ for every interesting index $i\in \{1,\dots,\ell |S|\}$. Indeed, if $x_i$ is a cousin of $y_i$, then it follows directly from the definition of $S'$ that $x_i\in S'$. If $x_i$ is not a cousin of $y_i$, then let $h(x_i)$ be the number of points $w\in S$ that are cousins of $x_i$. Since each of these points appears exactly $\ell$ times among $y_1,\dots,y_{\ell |S|}$, in this case there are $\ell\cdot h(x_i)$ indices $j\in \{1,\dots,\ell |S|\}\setminus \{i\}$ such that $y_j$ is a cousin of $x_i$. At most $\ell\cdot h(x_i)-1$ of these indices have the property that $x_j$ is a cousin of $x_i$, so there must be an index  $j\in \{1,\dots,\ell |S|\}\setminus \{i\}$ such that $y_j$ is a cousin of $x_i$ but $x_j$ is not a cousin of $x_i$. By transitivity of the cousin relation, this means that $x_j$ is not a cousin of $y_j$ and therefore the index $j$ cannot be boring. So $j$ is interesting and $x_i$ is a cousin of $y_j$, hence $x_i\in S'$.

So we have shown that $x_i\in S'$ for every interesting index $i\in \{1,\dots,\ell |S|\}$. Analogously, we can also show $z_i\in S'$ for every interesting index $i\in \{1,\dots,\ell |S|\}$.

As the set $S$ is mid-point-reducible, the non-empty subset $S'\su S$ must contain some non-mid-point $y\in S'$ of $S'$. By the definition of $S'$, for this point $y\in S$ there must be some interesting index $i\in \{1,\dots,\ell |S|\}$, such that $y$ is a cousin of $y_i$. As $y\in S'$ is a non-mid-point of $S'$, and $y_i\in S'$ is a cousin of $y$, the point $y_i$ is also a non-mid-point of $S'$. On the other hand, we have $2y_i=x_i+z_i$ and $x_i,z_i\in S'$, so by Definition \ref{defi-non-mid-point} this means that $x_i=z_i$. Hence $2y_i=x_i+z_i=2x_i$, and so $x_i=z_i$ is a cousin of $y_i$. So the index $i$ is boring, which is a contradiction. Thus, the set $A$ is indeed three-term progression free.
\end{proof}

\section{Construction of a non-mid-point reducible subset of \texorpdfstring{$\boldsymbol{\Z_m^2}$}{Zm2}}
\label{sect-constr-non-midpoint-deg-set}

In this section we prove Proposition \ref{proposition-2}, which is one of the main difficulties for proving our main results. Let us fix a positive integer $m$ throughout this section.

In order to construct the desired non-mid-point reducible subset $S\su \Z_m^2$, we map $\Z_m^2$ to a subset of the square $[0,1)^2$ in a certain way. We start by making some definitions and observations.

\begin{definition}
For real numbers $a,a'\in \mathbb{R}$, we write  $a\equiv a'\mod 1$ if $a-a'\in \mathbb{Z}$. For two points $(a,b),(a',b')\in \mathbb{R}^2$, we write $(a,b)\equiv (a',b')\mod 1$ if we have $a\equiv a'\mod 1$ and $b\equiv b'\mod 1$ (i.e. if we have $a-a'\in \mathbb{Z}$ and $b-b'\in \mathbb{Z}$).
\end{definition}

Note that for each point $(a',b')\in \mathbb{R}^2$, there exists exactly one point $(a,b)\in [0,1)^2$ with $(a,b)\equiv (a',b')\mod 1$.

\begin{definition}
For a finite subset $S\su [0,1)^2$, a point $y\in S$ is called a \emph{non-mid-point} of $S$ if there is no solution to the equation $2y\equiv x+z\mod 1$ with distinct $x,z\in S$. 
\end{definition}

Given real numbers $\alpha, \beta\in \mathbb{R}$, let us define a map $\varphi_{\alpha, \beta}:\Z_m^2\to [0,1)^2$ as follows. Every point in $\Z_m^2$ can be represented (uniquely) by a pair $(q,r)\in \{0,1,\dots,m-1\}^2$ of residue classes modulo $m$. Let us now define $\varphi_{\alpha, \beta}(q,r)=(a,b)$ for the unique point $(a,b)\in [0,1)^2$ with $(a,b)\equiv (\alpha+q/m,\beta+r/m)\mod 1$. Note that the map $\varphi_{\alpha, \beta}:\Z_m^2\to [0,1)^2$ is injective for any $\alpha, \beta\in \mathbb{R}$.

\begin{fact}\label{fact-equivalence}
Let $\alpha, \beta\in \mathbb{R}$ and $S\su \Z_m^2$, and consider the image $\varphi_{\alpha, \beta}(S)\su [0,1)^2$ of $S$ under the map $\varphi_{\alpha, \beta}:\Z_m^2\to [0,1)^2$. Then a point $y\in S$ is a non-mid-point of $S$ if and only if its image $\varphi_{\alpha, \beta}(y)$ is a non-mid-point of $\varphi_{\alpha, \beta}(S)$.
\end{fact}

\begin{proof}
Let $y\in S\su \Z_m^2$ be represented by the pair $(q_y,r_y)\in \{0,1,\dots,m-1\}^2$ of residue classes modulo $m$. Now, for any points $x,z\in S$, represented by pairs $(q_x,r_x),(q_z,r_z)\in \{0,1,\dots,m-1\}^2$ of residue classes, we have the chain of equivalences
\begin{align*}
2y=x+z\quad &\Leftrightarrow \quad 2q_y\equiv q_x+q_z\mod m\quad\text{ and }\quad 2r_y\equiv r_x+r_z\mod m\\
&\Leftrightarrow \quad 2q_y/m\equiv q_x/m+q_z/m\mod 1\quad\text{ and }\quad 2r_y/m\equiv r_x/m+r_z/m\mod 1\\
&\Leftrightarrow \quad 2(\alpha+q_y/m,\beta+r_y/m)\equiv (\alpha+q_x/m,\beta+r_x/m)+(\alpha+q_z/m, \beta+r_z/m)\mod 1\\
&\Leftrightarrow \quad 2\varphi_{\alpha, \beta}(y)\equiv \varphi_{\alpha, \beta}(x)+\varphi_{\alpha, \beta}(z)\mod 1.
\end{align*}
Thus, there is a solution to $2y=x+z$ with distinct $x,z\in S$ if and only if there is a solution to $\quad 2\varphi_{\alpha, \beta}(y)\equiv \varphi_{\alpha, \beta}(x)+\varphi_{\alpha, \beta}(z)\mod 1$ with distinct $x,z\in S$. Since $\varphi_{\alpha, \beta}$ is injective, this shows that $y$ is a non-mid-point of $S$ if and only if $\varphi_{\alpha, \beta}(y)$ is a non-mid-point of $\varphi_{\alpha, \beta}(S)$.
\end{proof}

The subset $T\su [0,1)^2$ defined in the following definition plays a crucial role in our proof. Using Fact \ref{fact-equivalence}, we will construct non-mid-point reducible subsets of $\Z_m^2$ from this subset $T\su [0,1)^2$.

\begin{definition}\label{defi-T}
Let us define
\[T_1=\Bigg\lbrace (a,b)\in \Big[0,\frac{1}{2}\Big)\times \Big[\frac{1}{2},1\Big)\,\Bigg|\,\frac{7}{12}\le a+b\le\frac{4}{3}\Bigg\rbrace \cup \Bigg\lbrace (a,b)\in \Big[0,\frac{1}{2}\Big)^2\,\Bigg|\,a+b>\frac{5}{6}\Bigg\rbrace\]
and
\[T_2=\Bigg\lbrace (a,b)\in \Big[\frac{1}{2},1\Big)\times \Big[0,\frac{1}{2}\Big)\,\Bigg|\,\frac{7}{12}\le a+b<\frac{5}{6} \text{ and }2a+b< \frac{3}{2}\Bigg\rbrace,\]
as well as $T=T_1\cup T_2\su [0,1)^2$.
\end{definition}

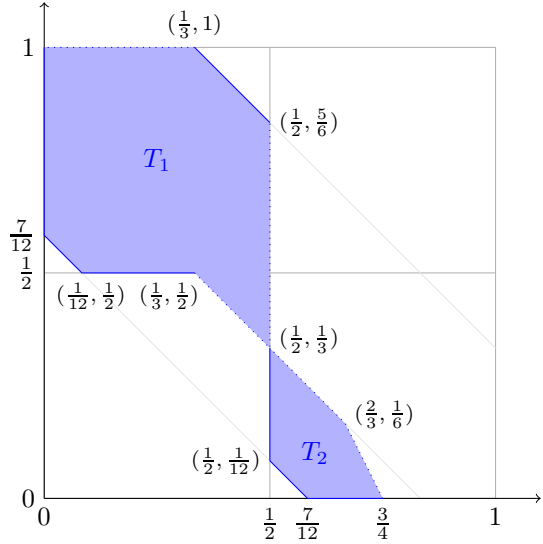
\begin{figure}
    \centering
\begin{tikzpicture}[scale=6]
\draw[gray!60] (0,0) -- (0,1) -- (1,1) -- (1,0) -- (0,0);
\draw[gray!60] (0,1/2) -- (1,1/2);
\draw[gray!60] (1/2,0) -- (1/2,1);

\draw[gray!20] (1/12,1/2) -- (1/2,1/12);
\draw[gray!20] (1/2,5/6) -- (1/1,1/3);
\draw[gray!20] (2/3,1/6) -- (5/6,0);

\draw[->] (0,0) -- (0,1.1);
\draw[->] (0,0) -- (1.1,0);
\draw (0,0) node[anchor=north]{0};
\draw (7/12,0) node[anchor=north]{$\frac{7}{12}$};
\draw (1/2,0) node[anchor=north]{$\frac{1}{2}$};
\draw (1,0) node[anchor=north]{1};
\draw (0,0) node[anchor=east]{0};
\draw (0,7/12) node[anchor=east]{$\frac{7}{12}$};
\draw (3/4,0) node[anchor=north]{$\frac{3}{4}$};
\draw (0,1/2) node[anchor=east]{$\frac{1}{2}$};
\draw (0,1) node[anchor=east]{1};
\draw (1/2,1/3+0.02) node[anchor=west]{\footnotesize{$(\frac{1}{2},\frac{1}{3})$}};
\draw (1/12+0.02,1/2) node[anchor=north]{\footnotesize{$(\frac{1}{12},\frac{1}{2})$}};
\draw (1/3+0.03,1/2) node[anchor=north east]{\footnotesize{$(\frac{1}{3},\frac{1}{2})$}};
\draw (1/2,1/12) node[anchor=east]{\footnotesize{$(\frac{1}{2},\frac{1}{12})$}};
\draw (2/3,1/6+0.02) node[anchor=west]{\footnotesize{$(\frac{2}{3},\frac{1}{6})$}};
\draw (1/2,5/6) node[anchor=west]{\footnotesize{$(\frac{1}{2},\frac{5}{6})$}};
\draw (1/3,1) node[anchor=south]{\footnotesize{$(\frac{1}{3},1)$}};
\fill[blue!30]  (0,1) -- (0,7/12) -- (1/12,1/2) -- (1/3,1/2) -- (1/2,1/3) -- (1/2,5/6) -- (1/3,1) -- (0,1);
\draw[blue]  (0,1) -- (0,7/12) -- (1/12,1/2) -- (1/3,1/2);
\draw[blue,dotted] (1/3,1/2) -- (1/2,1/3) -- (1/2,5/6);
\draw[blue]  (1/2,5/6) -- (1/3,1);
\draw[blue,dotted] (1/3,1) -- (0,1);
\fill[blue!30]  (1/2,1/3) -- (1/2,1/12) -- (7/12,0) -- (3/4,0) -- (2/3,1/6) --  (1/2,1/3);
\draw[blue] (1/2,1/3) -- (1/2,1/12) -- (7/12,0) -- (3/4,0);
\draw[blue,dotted] (3/4,0) -- (2/3,1/6) --  (1/2,1/3);
\draw[blue] (1/4,3/4) node {$T_1$};
\draw[blue] (6/10,1/10) node {$T_2$};
\end{tikzpicture}
    \caption{The set $T$ defined in Definition \ref{defi-T}}
    \label{figure-set-T}
\end{figure}

For an illustration of these sets, see Figure \ref{figure-set-T}. Note that the area of the set $T$ can be computed as follows: First, the area of $T_1$ is given by
\[\Big(\frac{1}{2}\Big)^2-\frac{1}{2}\cdot \Big(\frac{1}{12}\Big)^2-\frac{1}{2}\cdot \Big(\frac{1}{6}\Big)^2+\frac{1}{2}\cdot \Big(\frac{1}{6}\Big)^2=\Big(\frac{1}{2}\Big)^2-\frac{1}{2}\cdot \Big(\frac{1}{12}\Big)^2=\frac{72-1}{288}=\frac{71}{288}.\]
Next, the area of $T_2$ is given by
\[\frac{1}{2}\cdot \Big(\frac{1}{3}\Big)^2-\frac{1}{2}\cdot \Big(\frac{1}{12}\Big)^2-\frac{1}{2}\cdot \Big(\frac{1}{12}\Big)\cdot \Big(\frac{1}{6}\Big)=\frac{16-1-2}{288}=\frac{13}{288}.\]
Thus, the area of $T$ equals
\[\frac{71}{288}+\frac{13}{288}=\frac{84}{288}=\frac{7}{24}.\]

The key step for proving Proposition \ref{proposition-2} is to prove the following lemma concerning this subset $T\su [0,1)^2$.

\begin{lemma}\label{lemma-key}
Let $T\su [0,1)^2$ be defined as in Definition \ref{defi-T}. Then for any non-empty finite subset $S\su T$, there is a non-mid-point $y\in S$ of $S$.
\end{lemma}

Before proving this lemma, let us fist show how it implies Proposition \ref{proposition-2}. First, we obtain the following corollary  from the lemma.

\begin{corollary}\label{coro-key-lemma}
Let $T\su [0,1)^2$ be defined as in Definition \ref{defi-T}. Then for any real numbers $\alpha, \beta\in \mathbb{R}$, the preimage $\varphi_{\alpha, \beta}^{-1}(T)\su \Z_m^2$ of $T$ under the map $\varphi_{\alpha, \beta}:\Z_m^2\to [0,1)^2$ is a non-mid-point reducible subset of $\Z_m^2$.
\end{corollary}

\begin{proof}
For every non-empty subset $S'\su \varphi_{\alpha, \beta}^{-1}(T)$ the image $\varphi_{\alpha, \beta}(S')\su T$ is a non-empty finite subset of $T$. By Lemma \ref{lemma-key}, there exists a non-mid-point of  $\varphi_{\alpha, \beta}(S')$, so there is a point $y\in S'$ such that $\varphi_{\alpha, \beta}(y)$ is a non-mid-point of $\varphi_{\alpha, \beta}(S')$. By Fact \ref{fact-equivalence}, the point $y\in S'$ is then a non-mid-point of $S'$. So $\varphi_{\alpha, \beta}^{-1}(T)$ is non-mid-point reducible by Definition \ref{defi-non-mid-point-reducible}.
\end{proof}

The last missing piece to deduce Proposition \ref{proposition-2} is the following simple consequence of the fact that $T$ has area $7/24$.

\begin{fact}\label{fact-area-consequence}
Let $T\su [0,1)^2$ be defined as in Definition \ref{defi-T}. There exist real numbers $\alpha, \beta\in \mathbb{R}$, such that the preimage $\varphi_{\alpha, \beta}^{-1}(T)$ of $T$ under the map $\varphi_{\alpha, \beta}:\Z_m^2\to [0,1)^2$ has size $|\varphi_{\alpha, \beta}^{-1}(T)|\ge (7/24)\cdot m^2$.
\end{fact}

\begin{proof}
Let $\alpha, \beta\in [0,1]$ be independent uniform random variables in the interval $[0,1]$. Then for every point $x\in \Z_m^2$, the probability of having $\varphi_{\alpha, \beta}(x)\in T$ is precisely the area of $T\su [0,1)^2$, which is $7/24$. Hence
\[\mathbb{E}[|\varphi_{\alpha, \beta}^{-1}(T)|]=\sum_{x\in \Z_m^2}\mathbb{P}[x\in \varphi_{\alpha, \beta}^{-1}(T)]=\sum_{x\in \Z_m^2}\mathbb{P}[\varphi_{\alpha, \beta}(x)\in T]=\sum_{x\in \Z_m^2}7/24=(7/24)\cdot m^2,\]
which implies that for some choice of $\alpha, \beta\in [0,1]$ we must have $|\varphi_{\alpha, \beta}^{-1}(T)|\ge (7/24)\cdot m^2$.
\end{proof}

Now, combining Corollary \ref{coro-key-lemma} and Fact \ref{fact-area-consequence} immediately implies Proposition \ref{proposition-2}:

\begin{proof}[Proof of Proposition \ref{proposition-2}]
Let $T\su [0,1)^2$ be defined as in Definition \ref{defi-T}, and choose $\alpha, \beta\in \mathbb{R}$ as in Fact \ref{fact-area-consequence} such that $|\varphi_{\alpha, \beta}^{-1}(T)|\ge (7/24)\cdot m^2$. By Corollary \ref{coro-key-lemma}, the set $\varphi_{\alpha, \beta}^{-1}(T)\su \Z_m^2$ is non-mid-point reducible, so there exists a non-mid-point reducible subset of $\Z_m^2$ of size at least $(7/24)\cdot m^2$.
\end{proof}

It remains to prove Lemma \ref{lemma-key}. We first record the following simple facts that will be used in our proof.

\begin{fact}\label{fact-simple-properties-functions}
For the sets $T_1$, $T_2$ and $T=T_1\cup T_2$ defined in Definition \ref{defi-T}, the following statements hold:
\begin{itemize}
\item[(i)] For all $(a,b)\in T$, we have $7/12\le a+b\le 4/3$.
\item[(ii)] For all $(a,b)\in T_1$ and $(a',b')\in T_2$ with $a+b=a'+b'$, we have $a+a'< 1$.
\end{itemize} 
\end{fact}

\begin{proof}
Note that (i) follows immediately from the definition of $T$ (noting in particular that for all $(a,b)\in [0,1/2)^2$ with $a+b>5/6$ we have $7/12<5/6<a+b<1/2+1/2<4/3$).

For (ii), consider $(a,b)\in T_1$ and $(a',b')\in T_2$ with $a+b=a'+b'$. By the definition of $T_2$ we have $a'+b'< 5/6$, and hence $a+b< 5/6$. By the definition of $T_1$ this implies that $(a,b)\not\in [0,1/2)^2$ and consequently $b\ge 1/2$. Thus, we obtain
\[a=a'+b'-b\le a'+b'-\frac{1}{2}\]
and therefore
\[a+a'\le 2a'+b'-\frac{1}{2}<\frac{3}{2}-\frac{1}{2}=1,\]
where the last inequality follows from the definition of $T_2$.
\end{proof}

\begin{fact}\label{fact-simple-properties-midpoint}
Let $T\su [0,1)^2$ be defined as in Definition \ref{defi-T}. Consider points $(a,b), (a',b') \in T$, and let $(a^*,b^*)$  be their mid-point (given by $a^*=(a+a')/2$ and $b^*=(b+b')/2$). Then
\[a^*+b^*\le \frac{4}{3}.\]
Furthermore, if $a^*\ge 1/2$, then we have
\[a^*+b^*< \frac{13}{12}.\]
\end{fact}
\begin{proof}
For the first part, we simply observe that by Fact \ref{fact-simple-properties-functions}(i) we have
\[a^*+b^*=\frac{a+a'}{2}+\frac{b+b'}{2}=\frac{a+b}{2}+\frac{a'+b'}{2}\le \frac{4/3}{2}+\frac{4/3}{2}=\frac{4}{3}.\]
For the second part, let us assume that $a^*\ge 1/2$. Then we must have $a\ge 1/2$ or $a'\ge 1/2$, so let us assume $a'\ge 1/2$ without loss of generality. As $T_1\su [0,1/2)\times [0,1)$, this implies $(a',b')\in T_2$ and hence $a'+b'<5/6$ by the definition of $T_2$. Combining this with Fact \ref{fact-simple-properties-functions}(i), we obtain
\[a^*+b^*=\frac{a+a'}{2}+\frac{b+b'}{2}=\frac{a+b}{2}+\frac{a'+b'}{2}< \frac{4/3}{2}+\frac{5/6}{2}=\frac{2}{3}+\frac{5}{12}=\frac{13}{12}.\qedhere\]
\end{proof}

\begin{fact}\label{fact-AP-properties-functions}
Let $T\su [0,1)^2$ be defined as in Definition \ref{defi-T}. Consider points $ x,y,z\in T$ such that $2y\equiv x+z\mod 1$. Then, writing $x=(a_x,b_x)$, $y=(a_y,b_y)$ and $z=(a_z,b_z)$, we have
\begin{equation}
a_y+b_y\ge \frac{a_x+b_x}{2}+\frac{a_z+b_z}{2}.\label{eq-AP-property-sum}
\end{equation}
\end{fact}
\begin{proof}
Let $(a^*,b^*)$  be the mid-point of $x=(a_x,b_x)$ and $z=(a_z,b_z)$ (given by $a^*=(a_x+a_z)/2$ and $b^*=(b_x+b_z)/2$). Note that
\[(2a_y,2b_y)=2y\equiv x+z=(a_x,b_x)+(a_z,b_z)=(a_x+a_z,b_x+b_z)=(2a^*,2b^*)\mod 1,\]
meaning that $2a_y\equiv 2a^*\mod 1$ and $2b_y\equiv 2b^*\mod 1$. Hence $a^*-a_y\in \{-1/2,0,1/2\}$ and $b^*-b_y\in \{-1/2,0,1/2\}$, and consequently $(a^*+b^*)-(a_y+b_y)\in \{-1,-1/2,0,1/2,1\}$.

Furthermore, note that $a^*+b^*=(a_x+a_z)/2+(b_x+b_z)/2$ equals the right-hand side of (\ref{eq-AP-property-sum}), and so it suffices to show that $a_y+b_y \ge a^*+b^*$. So let us assume for contradiction that $a^*+b^*> a_y+b_y$. Then we automatically have $a^*+b^*\ge a_y+b_y+1/2$ (recalling that $(a^*+b^*)-(a_y+b_y)\in \{-1,-1/2,0,1/2,1\}$).

Since $b^*-b_y\le 1/2$, we must therefore have $a^*\ge a_y$. Similarly, we must have $b^*\ge b_y$, since $a^*-a_y\le 1/2$.

If $a^*\ge 1/2$, by the second part of Fact \ref{fact-simple-properties-midpoint} we have
\[a^*+b^*< \frac{13}{12}=\frac{7}{12}+\frac{1}{2}\le a_y+b_y+\frac{1}{2},\]
where the last inequality follows from Fact \ref{fact-simple-properties-functions}(i). This contradicts $a^*+b^*\ge a_y+b_y+1/2$, so we must have $a^*< 1/2$.

Now, we obtain $0\le a_y\le a^*<1/2$. Together with $a^*-a_y\in \{-1/2,0,1/2\}$, this implies $a^*=a_y$. So from $a^*+b^*\ge a_y+b_y+1/2$, we obtain $b_y\le b^*-1/2<1-1/2=1/2$.

Thus, we have $a_y<1/2$ and $b_y<1/2$, meaning that $(a_y,b_y)\in [0,1/2)^2$. By the definition of $T$, this implies $a_y+b_y>5/6$ and hence
\[a^*+b^*\le \frac{4}{3}=\frac{5}{6}+\frac{1}{2}< a_y+b_y+\frac{1}{2},\]
where the first inequality follows from the first part of Fact \ref{fact-simple-properties-midpoint}. This is again a contradiction to $a^*+b^*\ge a_y+b_y+1/2$. 
\end{proof}

For our last fact, we need some more notation. Let us define the function $g:[0,1)\to [0,1/2)$ by setting
\[g(t)=\begin{cases}
t &\text{if }t\in [0,1/2)\\
t-1/2 &\text{if }t\in [1/2,1)
\end{cases}\]
for all $t\in [0,1)$. Note that for $t,t'\in [0,1)$ with $2t\equiv 2t'\mod 1$, we always have $g(t)=g(t')$.

\begin{fact}\label{fact-AP-properties-function-g}
Let $T\su [0,1)^2$ be defined as in Definition \ref{defi-T}. Consider points $ x,y,z\in T$ such that $2y\equiv x+z\mod 1$. Let us write $x=(a_x,b_x)$, $y=(a_y,b_y)$ and $z=(a_z,b_z)$, and assume that $a_x+b_x=a_z+b_z$. Then, we have
\[g(a_y)\ge \frac{g(a_x)+g(a_z)}{2}.\]
Furthermore, if $g(a_x)=g(a_y)=g(a_z)$, then $x=z$.
\end{fact}
\begin{proof}
As in the previous proof, let $(a^*,b^*)$  be the mid-point of $x=(a_x,b_x)$ and $z=(a_z,b_z)$ (given by $a^*=(a_x+a_z)/2$ and $b^*=(b_x+b_z)/2$) and note that $2a_y\equiv 2a^*\mod 1$ and $2b_y\equiv 2b^*\mod 1$. Hence $g(a_y)=g(a^*)$.

Recall that $x,z\in T=T_1\cup T_2$. We now distinguish three cases, first the case that $x,z\in T_1$, second the case that $x,z\in T_2$, and third the case that $x\in T_1$ and $z\in T_2$. Note that the fourth possible case that $x\in T_2$ and $z\in T_1$ is analogous to the third case (since we can interchange the roles of $x$ and $z$ in the statement).

\textit{Case 1: $x,z\in T_1$.} In this case, we have $a_x,a_z\in [0,1/2)$ and hence also $a^*=(a_x+a_z)/2\in [0,1/2)$. Therefore $g(a_x)=a_x$, $g(a_z)=a_z$ and $g(a^*)=a^*$, and consequently
\[g(a_y)=g(a^*)=a^*=\frac{a_x+a_z}{2}=\frac{g(a_x)+g(a_z)}{2}.\]
If furthermore $g(a_x)=g(a_y)=g(a_z)$, then $a_x=g(a_x)=g(a_z)=a_z$. Together with $a_x+b_x=a_z+b_z$, this also implies $b_x=b_z$, so $x=(a_x,b_x)=(a_z,b_z)=z$.

\textit{Case 2: $x,z\in T_2$.} In this case, we have $a_x,a_z\in [1/2,1)$ and hence also $a^*=(a_x+a_z)/2\in [1/2,1)$. Therefore $g(a_x)=a_x-1/2$, $g(a_z)=a_z-1/2$ and $g(a^*)=a^*-1/2$, and consequently
\[g(a_y)=g(a^*)=a^*-1/2=\frac{a_x+a_z-1}{2}=\frac{g(a_x)+g(a_z)}{2}.\]
If furthermore $g(a_x)=g(a_y)=g(a_z)$, then $a_x=g(a_x)+1/2=g(a_z)+1/2=a_z$. Together with $a_x+b_x=a_z+b_z$, this also implies $b_x=b_z$, so $x=(a_x,b_x)=(a_z,b_z)=z$.

\textit{Case 3: $x\in T_1$ and $z\in T_2$.} In this case, we have $a_x\in [0,1/2)$ and $a_z\in [1/2,1)$. Furthermore, by Fact \ref{fact-simple-properties-functions}(ii) we have $a_x+a_z<1$ and hence $a^*=(a_x+a_z)/2<1/2$. Therefore $g(a_x)=a_x$ and $g(a^*)=a^*$, while $g(a_z)=a_z-1/2$. Thus we obtain
\[g(a_y)=g(a^*)=a^*=\frac{a_x+a_z}{2}=\frac{g(a_x)+g(a_z)+1/2}{2}>.\frac{g(a_x)+g(a_z)}{2},\]
implying the desired inequality and showing that $g(a_x)=g(a_y)=g(a_z)$ is not possible in this case (so the second part of the fact is tautologically true in this case).
\end{proof}

Finally, let us prove Lemma \ref{lemma-key}.

\begin{proof}[Proof of Lemma \ref{lemma-key}]
Let $S\su T$ be a non-empty finite subset of $T$. Let $\mu=\min_{(a,b)\in S} a+b$ be the minimal value of $a+b$ for all points $(a,b)\in S$. Furthermore, define $\nu=\min \{g(a)\mid (a,b)\in S, a+b=\mu\}$ to be the minimal value of $g(a)$ among all points $(a,b)\in S$ with $a+b=\mu$. Finally, let us choose a point $y=(a_y,b_y)\in S$ with $a_y+b_y=\mu$ and $g(a_y)=\nu$. Our goal is to show that $y$ is a non-mid-point of $S$.

Suppose for contradiction that $y$ is not a non-mid-point of $S$, then there exist distinct points $x=(a_x,b_x)\in S$ and $z=(a_z,b_z)\in S$ with $2y\equiv x+z\mod 1$. Now, by Fact \ref{fact-AP-properties-functions} we have
\[\mu=a_y+b_y \ge \frac{a_x+b_x}{2}+\frac{a_z+b_z}{2},\]
but on the other hand by the definition of $\mu$ we must have $a_x+b_x\ge \mu$ and $a_z+b_z\ge \mu$. Thus, we obtain $a_x+b_x=a_y+b_y=a_z+b_z= \mu$. By Fact \ref{fact-AP-properties-function-g}, this implies that 
\[\nu=g(a_y)\ge \frac{g(a_x)+g(a_z)}{2},\]
but on the other hand by the definition of $\nu$ we must have $g(a_x)\ge \nu$ and $g(a_z)\ge \nu$ (recalling that $a_x+b_x= \mu$ and $a_z+b_z= \mu$). Consequently, we obtain $g(a_x)=g(a_y)=g(a_z)= \nu$. Now, the second part of Fact \ref{fact-AP-properties-function-g} shows that $x=z$, which is a contradiction. This contradiction shows that $y\in S$ is indeed a non-mid-point of $S$.
\end{proof}

\section{Concluding remarks}
\label{sect-concluding}

Understanding the size of the largest three-term progression free subset of $\Fpn$ for a fixed prime $p\ge 3$ and large $n$ is still an open problem. Ellenberg and Gijswijt \cite{ellenberg-gijswijt} proved an upper bound of $(c_p p)^n$ for some constant $c_p<1$ depending on $p$ (with $c_p\to 0.841...$ for $p\to \infty$). Our results show a lower bound of the form $(\sqrt{7/24}\cdot p)^{n-O(\log n)}$ for any fixed prime $p\ge 3$. It is plausible that with a similar approach one can get a better constant than $\sqrt{7/24}\approx 0.54006...$, potentially by applying Proposition \ref{proposition-1} for some $d>2$. Of course, the crucial difficulty with this approach is finding a non-mid-point reducible subset of $\Fp^d$ of size larger than $(\sqrt{7/24}\cdot p)^{d}$ for some $d$.

In order to find such a set, it could conceivably be helpful to relax the reducibility condition for the desired set $S\su \Fp^d$ (see Definition \ref{defi-non-mid-point-reducible}). Indeed, instead of demanding that every non-empty subset $S'\su S$ contains a non-mid-point (i.e.\ a point that is not the middle term of an arithmetic progression in $S'$), it is sufficient to demand that every non-empty subset $S'\su S$ contains a point that is not the middle term of an arithmetic progression in $S'$ or a point that is not the first or third term of an arithmetic progression in $S'$. See \cite[Proof of Theorem 3.13]{elsholtz-pach} and \cite[Proof of Theorem 2.3]{elsholtz-klahn-lipnik}, where such arguments were used (and see also \cite{elsholtz-lipnik}, where in a related situation a further extension based on integer programming was used).

For a fixed choice of a (small) prime $p\ge 3$, or more generally for a fixed choice of a (small) integer $m\ge 3$, it may be possible to obtain a non-mid-point reducible subset $S\su \Z_m^2$ of size more than $(7/24)\cdot m^2$ (i.e.\ larger than the size guaranteed in Proposition \ref{proposition-2}), by applying Corollary \ref{coro-key-lemma} to the set $S=\varphi_{\alpha, \beta}^{-1}(T)\su \Z_m^2$ for some smartly chosen values of $\alpha, \beta\in \mathbb{R}$ (rather than taking $\alpha$ and $\beta$ randomly as in the proof of Fact \ref{fact-area-consequence}). However, for larger $m$, the size of the set $\varphi_{\alpha, \beta}^{-1}(T)$ is close to $(7/24)\cdot m^2$ for any choice of $\alpha, \beta\in \mathbb{R}$. In particular, for sufficiently large $m$, one can make the construction of the set $A$ in Theorems \ref{thm-main} and \ref{thm-main-2} completely explicit by taking $S=\varphi_{0, 0}^{-1}(T)$ in the proof of Proposition \ref{proposition-2}.

It is well-known that for every fixed prime $p\ge 3$, the maximum size of a three-term progression free subset of $\Fpn$ is of the form $(\gamma_p p)^{n-o(n)}$ for all $n$, where $\gamma_p$ is a constant only depending on $p$ (see e.g. \cite[Proposition 3.8]{elsholtz-pach}). Ellenberg and Gijswijt \cite{ellenberg-gijswijt} showed that $\gamma_p<1$ for all $p$ and that $\limsup_{p\to \infty}\gamma_p<0.842$ (see \cite[Eq.\ (4.11)]{blasiak-et-al}). Our work shows that $\gamma_p>0.54$ for all $p$ and hence in particular that $\liminf_{p\to \infty}\gamma_p>1/2$. It is an open question whether the limit $\lim_{p\to \infty}\gamma_p$ exists, see \cite[Conjecture 3.9]{elsholtz-pach}.

Finally, we remark that the known constructions in the integer case, i.e.\ the constructions of large subsets of $\{0,\dots,N\}$ without three-term arithmetic progressions, are also of a high-dimensional nature. In these constructions, one chooses integers $t$ and $n$ such that $(2t-1)^n\le N$. One can then obtain a three-term progression free subset of $\{0,\dots,N\}$ from a three-term progression free subset of $\{0,\dots,t-1\}^n\su \Z^n$ by considering $n$-tuples in $\{0,\dots,t-1\}^n\su \Z^n$ as numbers written in base $2t-1$ (for three such numbers $x,y,z$ there cannot be any carries in the equation $2y=x+z$ written in base $2t-1$, and hence any three-term arithmetic progression of such numbers would correspond to a three-term progression in the original subset of $\{0,\dots,t-1\}^n$). The earliest such construction, described by Erd\H{o}s and Tur\'{a}n \cite{erdos-turan} and attributed by them to Szekeres, used this idea for $t=2$, taking the three-term progression free set $\{0,1\}^n$ (the resulting subset of $\{0,\dots,N\}$ then consists of all integers whose base-3 representation does not contain the digit $2$). Szekeres  also conjectured that this construction is best possible, which turned out to be incorrect. Salem and Spencer \cite{salem-spencer} generalized Szekeres' idea, by considering subsets of $\{0,\dots,t-1\}^n$ with those points where each of the numbers  $0,\dots,t-1$ occurs in a prescribed number of coordinates, i.e.~when written as integers in base $m=2t-1$, each of the digits $0,\dots, t-1$ occurs a prescribed number of times. Finally, Behrend's construction \cite{behrend} is based on taking the intersection of $\{0,\dots,t-1\}^n$ with a suitably chosen sphere (by convexity of the sphere, such an intersection cannot contain any three-term arithmetic progressions). Elkin \cite{elkin} slightly improved Behrend's lower bound for the size of the largest three-term progression free subset of $\{0,\dots,N\}$ to a bound on the order of $N(\log N)^{1/4}2^{-2\sqrt{2\log_2 N}}$ by taking a subset of $\{0,\dots,t-1\}^n$ consisting of points very close to a given sphere. 

Given the high-dimensional nature of all these constructions for progression-free sets in the integers, relying on three-term progression free subsets of $\{0,\dots,t-1\}^n\su \Z^n$, it is natural to wonder whether the ideas in this paper might also be useful for this integer version of the problem. However, one crucial difficulty is that in order to avoid carries when performing the additions in the equation $2y=x+z$, one seems to need to restrict oneself to only half of the available digits in base $2t-1$.

MSC classification:\\
05D05 Extremal set theory\\
11B25 Arithmetic progressions


\begin{thebibliography}{99}

\bibitem{alon-shpilka-umans} N. Alon, A. Shpilka, and C. Umans, \textit{On sunflowers and matrix multiplication}, Comput. Complexity \textbf{22}
(2013), 219--243.

\bibitem{behrend} F. A. Behrend, \textit{On sets of integers which contain no three terms in arithmetical progression}, Proc. Nat. Acad. Sci. U.S.A. \textbf{32} (1946), 331--332.

\bibitem{blasiak-et-al} J. Blasiak, T. Church, H. Cohn, J. A. Grochow, E. Naslund, W. F. Sawin, and C. Umans, \textit{On cap sets
and the group-theoretic approach to matrix multiplication}, Discrete Anal. 2017, Paper No. 3, 27pp.

\bibitem{bloom-sisask} T. F. Bloom and O. Sisask, \textit{An improvement to the Kelley-Meka bounds on three-term arithmetic progressions}, preprint, 2023, arXiv:2309.02353.

\bibitem{brown-buhler} T. C. Brown and J. P. Buhler, \textit{A density version of a geometric Ramsey theorem}, J. Combin. Theory Ser. A \textbf{32} (1982), 20--34.

\bibitem{calderbank-fishburn} A. R. Calderbank and P. C. Fishburn, \textit{Maximal three-independent subsets of $\{0,1,2\}^n$}, Des. Codes Cryptogr. \textbf{4} (1994), 203--211.

\bibitem{croot-lev-pach} E. Croot, V. F. Lev, and P. P. Pach, \textit{Progression-free sets in $\Z_4^n$ are exponentially small}, Ann. of Math. \textbf{185} (2017), 331--337.

\bibitem{edel} Y. Edel, \textit{Extensions of generalized product caps}, Des. Codes Cryptogr. \textbf{31} (2004), 5--14.

\bibitem{elkin} M. Elkin, \textit{An improved construction of progression-free sets}, Israel J. Math. \textbf{184} (2011), 93--128.

\bibitem{ellenberg-gijswijt} J. S. Ellenberg and D. Gijswijt, \textit{On large subsets of $\mathbb{F}_q^n$ with no three-term arithmetic progression}, Ann. of Math. \textbf{185} (2017), 339--343.

\bibitem{elsholtz-et-al-line-free} C. Elsholtz, J. F\"uhrer, E. F\"uredi, B. Kov\'acs, P. P. Pach, D. G. Simon, and N. Velich, \textit{Maximal line-free sets in $\mathbb{F}_p^n$}, preprint, 2023, arXiv:2310.03382.


\bibitem{elsholtz-klahn-lipnik}
C. Elsholtz, B. Klahn, and G. F. Lipnik, \textit{Large subsets of $\mathbb{Z}_m^n$ without arithmetic progressions}, Des. Codes Cryptogr. \textbf{91} (2023), 1443--1452.

\bibitem{elsholtz-lipnik} C. Elsholtz and G. F. Lipnik, \textit{Exponentially larger affine and projective caps},
Mathematika \textbf{69} (2023), 232--249.

\bibitem{elsholtz-lipnik-siebenhofer}
C. Elsholtz, G. F. Lipnik, and M. Siebenhofer, \textit{New constructions of large caps}, chapter 5 in Ph.D. Thesis of G. F. Lipnik, On Integer Partitions, Caps, Progression-Free Sets and $q$-Regular Sequences, Graz University of Technology, March 2023.

\bibitem{elsholtz-pach} C. Elsholtz and P. P. Pach, \textit{Caps and progression-free sets in $\mathbb{Z}_m^n$}, Des. Codes Cryptogr. \textbf{88} (2020), 2133--2170.

\bibitem{erdos-turan} P. Erd\H{o}s and P. Tur\'an, \textit{On Some Sequences of Integers}, J. London Math. Soc. \textbf{11} (1936), 261--264.

\bibitem{fox-pham} J. Fox and H. T. Pham, \textit{Popular progression differences in vector spaces II}, Discrete Anal. 2019, Paper No. 16, 39 pp.

\bibitem{frankl-graham-rodl} P. Frankl, R. L. Graham, and V. Rödl, \textit{On subsets of abelian groups with no 3-term arithmetic progression}, J. Combin. Theory Ser. A \textbf{45} (1987), 157--161.

\bibitem{fuehrer}
J.~F\"uhrer, \textit{More on maximal line-free sets in $\mathbb{F}_p^n$}, submitted, 2023.

\bibitem{kelley-meka} Z. Kelley and R. Meka, \textit{Strong Bounds for 3-Progressions}, preprint, 2023, arXiv:2302.05537.

\bibitem{lev:2004}
V. F.~Lev,
\textit{Progression-free sets in finite abelian groups},
J. Number Theory \textbf{104} (2004), 162--169.


\bibitem{lin-wolf}
Y.~Lin and J.~Wolf, \textit{On subsets of $\mathbb{F}_q^n$ containing no $k$-term progressions}, European J. Combin. \textbf{31} (2010), 1398--1403.

\bibitem{meshulam}
R. Meshulam,
\textit{On subsets of finite abelian groups with no 3-term arithmetic progressions},
J. Combin. Theory Ser. A \textbf{71} (1995), 168--172.

\bibitem{naslund}
E.~Naslund, \textit{Lower bounds for the Shannon Capacity of Hypergraphs}, manuscript, 2024.


\bibitem{peluse} S. Peluse, \textit{Finite field models in arithmetic combinatorics -- twenty years on}, preprint, 2023, arXiv:2312.08100.

\bibitem{google-deep-mind} B. Romera-Paredes, M. Barekatain, A. Novikov, M. Balog, M. P. Kumar, E. Dupont, F. J. R. Ruiz, J.~S.~Ellenberg, P. Wang, O. Fawzi, P. Kohli, and A. Fawzi, \textit{Mathematical discoveries from program search with large language models}, Nature, to appear.

\bibitem{salem-spencer} R. Salem and D. C. Spencer, \textit{On sets of integers which contain no three terms in arithmetical progression},
Proc. Nat. Acad. Sci. U.S.A. \textbf{28} (1942), 561--563.




\bibitem{tao-vu} T. Tao and V. H . Vu, \textit{Additive Combinatorics}, Cambridge University Press, 2006.

\bibitem{tyrrell} F. Tyrrell, \textit{New lower bounds for cap sets}, Discrete Anal. 2023, Paper No. 20, 18pp.
\end{thebibliography}
\end{document}